\documentstyle[amssymb]{article}

\begin{document}
\vskip 10 pt

\noindent {\large \bf On $(p,q;\alpha, \beta,l)$-deformed
oscillators and their oscillator algebras}

\vskip 10 pt

{I. M. Burban}

Institute for Theoretical Physics, Kiev 03143, Ukraine

\vskip  40 pt

\begin{abstract}
We present a description of a new kind of the deformed canonical
commutation relations, their representations and generated by them
Heisenberg--Weyl algebra. This deformed algebra allows us to
define operations of the Hopf algebra structure: comultication,
counit and antipode. We discuss properties of a discrete spectrum of
the Hamiltonian of the deformed harmonic oscillator corresponding
to this oscillator-like system.
\end {abstract}

\vskip 35 pt

\noindent {\bf 1. Introduction}
\medskip

\noindent Attempts to deform the canonical commutation relations
are repeatedly undertaken in physical theories \cite {ACL}.

With the emergence of the quantum groups (quantum multi-parameter
deformed universal enveloping algebras of Lie algebras) became
evident their important role for theoretical and mathematical
physics. From the physical point of view, the interest to quantum
deformations, in particular to the quantum $(p,q)$-deformations of
Lie  algebras, is connected with the possible applications them in
the quantum field theory (the conformal, topological field
theories, etc.). The description of the two-parameter quantum
groups and their representations has been started in the works
\cite{CZ}, \cite{ShWZ}, \cite{SW}.

As in the classical case, the problem of realization of
$q$-deformed algebras by the one-parameter deformed creation and
annihilation operators (the Jordan-Schwinger construction)
\cite{B}, \cite {M} is important for representation theory of
quantum groups.

This problem remains important for $(p,q)$-deformed cases as well.
The further exploration of these deformations led to investigation
of $(p,q)$-deformed canonical commutation relations \cite {ChJ},
\cite {J}, \cite {Q}.

At the study of the quantum groups and algebras became evident
their connection with the noncommutative geometry, special
functions of $q$-analysis and others branches of mathematics.

In the framework of this program the problems of q-analysis and
q-special functions have got natural extension to the
$(p,q)$-case. Already in the paper \cite{ChJ} and in further works
\cite {BK}, \cite{JR} the $(p,q)$-exponential and other
(p,q)-deformed functions have been introduced and their properties
have been investigated.

In the paper \cite {BK} the definition of the basic $(p,q)$-
hypergeometric series ${}_r\Psi_{r-1}$ was given and  their
properties were investigated. In \cite{JR} the general
$(p,q)$-hypergeometric series are defined and various
$q$-identities are converted into their $(p,q)$-analogs. In this
framework the $(p,q)$-differentiation, and the $(p,q)$-Jackson
integration has been defined and their main properties has been
studied \cite {Br}, \cite{BK}.

 The problem of extension of these results to the generalized
one- and two-parameter deformed cases naturally arise. An example
of such generalized $q$-deformed algebra with the Hopf algebra
structure has been studied in \cite{C-C-N-U}, \cite{OS}.

In this paper we define a new kind of the deformed canonical
commutation relations and connected with them Heisenberg-Weyl
algebra. We study representations of this algebra. This
generalized deformed algebra allows us to define the operations of
the comultication, antipode and counit which satisfy the axioms of
the Hopf algebra structure. This generalized
$(p,q;\alpha,\beta,l)$-deformed system includes as a particular
case the systems of \cite{C-C-N-U}, \cite{OS}. We discuss
properties of discrete spectrum of the Hamiltonian of the deformed
harmonic oscillator corresponding to this oscillator-like system.

\medskip
\medskip
\medskip
{\bf 2. The $(p,q;\alpha,\beta,l)$-deformed oscillator algebra and
its representations}
\medskip
\medskip
\medskip

A deformed Heisenberg--Weyl algebra is defined as the associative
algebras generated by the operators $\{{\bf 1}, a, a^+, N \}$ and
defining relations
\begin{equation}\label{burban: basic1} [N, a]= - a,\quad [N, a^+]= a^+, \end{equation}
\begin{equation}\label{burban:basic1a}
a^+ a = f(N),\quad a a^+ = f(N+1),
\end{equation}
where structure function $f(x)$ is a positive analytic function.
Instead of (\ref{burban:basic1a}) one can considere the relation
\begin{equation}\label{burban: basic2} [a, a^+] = f(N+1) - \,f(N),\end{equation}
although the algebras in this two cases in general are not
isomorphic. We define the generalize deformed Heisenberg--Weyl
algebra as an associative algebra generated by generators $ {\bf
1}, a, a^+, N $ satisfying the defining relations
\begin{equation}\label{burban: basic1b} [N, a]= - l a,\quad [N, a^+]= l a^+, \end{equation}
\begin{equation}\label{burban: basic2b}
[a, a^+]_{A} = f(N+l) -{A}\,f(N),
\end{equation}
where $[a, a^+]_{A} = aa^+ -{A}\, a^+ a,$ and $A,l\in {\mathbb
R}.$

The structure functions $f(x)$ in (\ref{burban: basic2}) and
(\ref{burban: basic2b}) characterize the deformation scheme. For
various known deformations of the harmonic oscillator they are
given: $$ f(n) = \frac{1}{2}n \quad {\rm for\,\,the\,\,
oscillator\,\, of\,\,standard \, quantum\,}$$ mechanics; $$
f(n)=[n]\quad {\rm and}\quad f(n) = q^{\alpha n+ \beta}[n],\,{\rm
where}\quad[n] = \frac {1 - q^n}{1-q},$$ {\rm define the
Arik--Coon and its generalization;}

$$ f(n)= [n]\quad{\rm and}\quad f(n)=[\alpha n + \beta], \,{\rm
where}\quad{\rm and}\quad{\rm and}\quad [n]= \frac {q^{-n}-
q^{n}}{q^{-1}- q},$$ define the Biedenharn--Makfarlane and its
symmetric generalization;
$$ f(n) = [n]\quad{\rm and}\quad f(n )= [\alpha n + \beta],\, {\rm
where}\quad [n]= \frac{p^{-n}-q^n}{p^{-l}-q^l },$$ define the
two-parameter deformation and its symmetric generalization.

In comparison with the one-parameter deformed commutation
relations the multi-parameter deformation are less understood
\cite{CZ}, \cite{ShWZ}, \cite{SW}.

Two-parameter analogs of the one-parameter symmetric deformation
\cite{ChJ}, \cite {Q} of the oscillator algebra are defined as an
associative algebra generated by the operators ${\bf 1}, a, a^+, N
$ and defining relations
\begin{equation} aa^+ - q a^+ a = p^{-N},\end{equation}
\begin{equation} aa^+ - p^{-1}a^+ a = q ^{N},\end{equation}
\begin{equation} [N, a] = - a,\quad [N, a^+] =
a^+.\end{equation} The generalized the Biedenharn--Macfarlane
$q$-oscillator algebra with defining relations
\begin{equation} aa^+ - q a^+ a = q^{-\alpha N - \beta}, \end{equation}
\begin{equation}
aa^+ - q^{-1}a^+ a = q ^{\alpha N + \beta},
\end{equation}
\begin{equation}
[N, a] = - a,\quad [N, a^+] = a^+
\end{equation}
and its Hopf algebra structure have been studied in the papers
\cite{C-C-N-U} and \cite{OS}. The properties of this algebra and
of the corresponding deformed oscillator it were studied in \cite
{BDY}.

By analogy with the deformation of \cite{C-C-N-U} and \cite {OS}
we introduce {\it the corresponding
$(p,q;\alpha,\beta,l)$-deformed canonical commutation
relations.\rm} The \, $(p,q;\alpha, \beta,l)$-deformed oscillator
algebra is given by the generators ${\bf 1}, a, a^+, N $ and the
commutation relations
\begin{equation}\label{burban: rel1} aa^+ - q^l a^+ a = p^{-\alpha N-\beta},\end{equation}
\begin{equation}\label{burban: rel2} aa^+ - p^{-l}a^+ a = q
^{\alpha N + \beta},\end{equation}
\begin{equation}\label{burban: rel3a}[N, a] = -l\,a, \quad[N, a^+] = l
a^+,\end{equation} where the function $f(n)$ has the form
\begin{equation}
f(n)= \left(\frac{p^{-\alpha n - \beta }- q^{\alpha n +
\beta}}{p^{-l}-q^l} \right)
\end{equation}
with $\alpha,\,\beta,\,l\in {\mathbb R}.$

Instead of the relations (\ref{burban: rel1}) and (\ref{burban:
rel2}) we can consider the relations
\begin{equation}\label{burban: 14a} aa^+ = \frac{p^{-\alpha N - \beta -l}- q^{\alpha N +
\beta+l}}{p^{-l}- q^l},\quad a^+a = \frac{p^{-\alpha N + \beta}-
q^{\alpha N + \beta}}{p^{-l} - q^l},
\end{equation}
which together with relations (\ref{burban: rel3a}) define an
oscillator algebra which in general is not isomorphic to the
algebra defined above. The difficulties to supply it with a Hopf
algebra structure are the same as in \cite{OS}.

Nevertheless, if we will replace the relations (\ref{burban: 14a})
by
\begin{equation}\label{burban: rel3}
[a, a^+]_{A} = \frac{p^{- \alpha N - \beta -l}-  q^{\alpha N +
\beta + l}} {p^{-l} - q^l}- {A}\, \frac{p^{-\alpha N - \beta}-
q^{\alpha N + \beta}}{p^{-l}- q^l},
\end{equation}
then we obtain a new algebra which can be considered as
$(p,q;\alpha,\beta,l)$-deformed Heisenberg-Weyl algebra. This
algebra generated relations (\ref{burban: rel3a}) and
(\ref{burban: rel3}), as we shall show in next section, admits a
Hopf algebra structure for a properly chosen constant  $A.$

The representation of the creation and annihilation operators $a,
a^+$ and the operator of number particles $N$ of the relations
(\ref{burban: rel1}), (\ref{burban: rel2}), (\ref{burban: rel3a})
in the Hilbert space ${\cal H}$ with the basis $\{|n\rangle\}, n =
0,1,2 \ldots$ are defined as follows
\begin{equation} a|n\rangle =
\left(\frac{p^{-\alpha -\beta}-q^{\alpha +
\beta}}{p^{-l}-q^l}\right)^{1/2}|n-l\rangle,\quad a^+\,|n\rangle
=\left(\frac{p^{-\alpha -\beta -l}-q^{\alpha + \beta + l
}}{p^{-l}-q^l}\right)^{1/2}|n+l\rangle,
\end{equation}
\begin{equation}N\,|n\rangle = n\,|n\rangle.
\end{equation}

In the space of functions (analytic if $l/\alpha$ is integer
number) we can define the difference derivative
\begin {equation}
D f(z)= \frac{f(p^{-\alpha}z)p^{- \beta} -
f(q^{\alpha}z)q^{\beta}}{(p^{-l}- q^l)z^{l/\alpha}}.\end{equation}
It follows $$ Dz^n = \frac{p^{-\alpha n - \beta}- q^{\alpha n +
\beta}}{p^{-l}- q^l}\,z^{n-l/\alpha} = \frac{z^n}{z^{l/\alpha}}\,
\frac{p^{- \alpha n - \beta}- q^{\alpha n + \beta}}{p^{-l} -
q^l}\,\frac{1}{(n)!}\frac{d^{n}z^{n}}{d z^{n}} $$ and (if
$l/\alpha$ is integer number)
\begin{equation} D f(z) = \sum_{n=0}^{\infty}
a_n Dz^n = \sum_{n=1}^{\infty}\frac{z^n}
{z^{l/\alpha}}\frac{p^{-\alpha n-\beta} - q^{\alpha n +
\beta}}{p^{-l} -
q^l}\frac{1}{n!}\frac{d^n}{dz^n}f(z)\end{equation} for an analytic
function $ f(z)= \sum_{n=0}^{\infty}a_n z^n. $
\medskip

Now we can give a realization of the relations (\ref{burban:
 rel1}), (\ref{burban: rel2}), (\ref{burban: rel3a}) in this space by the operators
\begin{equation}\label{burban: an1} a: f\to D f, \end{equation}
\begin{equation}\label{burban: an2} a^+ :f\to z^{l/\alpha}f, \end{equation}
\begin{equation}\label{burban: an3} N: f\to \alpha\,z\frac{d}{d z}, \end{equation}
\begin{equation}\label{burban: an4}q^N: f\to q^{z\frac{d}{d z}}f = f(q z),\end{equation}
\begin{equation}\label{burban: an5} p^{-N}: f\to p^{-z\frac{d}{d z}}f = f(p^{-1}z).
\end{equation}
Indeed, from (\ref{burban: an1}) and (\ref{burban: an3}) we obtain
$$ N a^+ f(z) = \alpha z \frac{d}{d z}(z^{l/\alpha} f(z)) = l
z^{l/\alpha}f + \alpha z^{1 + l/\alpha}\frac{d}{d z} f(z)$$ and
$$a^+ N f(z) = \alpha z^{1+l/\alpha }\frac{d}{d z}f(z). $$ It
follows that \begin{equation} [N, a^+]f = l\,a^+ f.\end{equation}
Analogously, from (\ref{burban: an2}) and (\ref{burban: an3}) we
get
$$ N a f = - l\frac{ \left(
f(p^{-\alpha}z)p^{-\beta}-f(q^{\alpha}z)q^{\beta})\right)}{z^{l/\alpha}(p^{-l}-q^l)
} + \alpha \frac{z p^{-\alpha-\beta}f'(p^{-\alpha}z)- z q^{\alpha
+ \beta} f'(q^{\alpha }z)}{z^{l/\alpha}(p^{-l}-q^l)}$$ and $$ a N
f = \alpha \frac{z p^{-\alpha-\beta}f'(p^{-\alpha}z)- z q^{\alpha
+ \beta} f'(q^{\alpha }z)}{z^{l/\alpha}(p^{-l}-q^l)}.$$ It follows
that
\begin{equation}[N, a]= -l\,a. \end{equation}

In a similar way, from (\ref{burban: an1}), (\ref{burban: an2}) we
have $$ a^+ \,a f(z) = \frac{f(p^{-1}z)p^{-\beta}- f(q^{\beta}
z)q^{\beta}}{p^{-l} - q^l}$$ and $$ a\,a^+ f(z) =
\frac{f(p^{-\alpha}z)p^{-l-\beta}- f(q^{\alpha} z)q^{l +
\beta}}{p^{-l}-q^l}.$$ Therefore,
\begin{equation}a\,a^+ - q^l a^+\,a = p^{-\alpha N -\beta},
\quad a\,a^+ -p^{-l} a^+\,a = q^{\alpha N + \beta}. \end{equation}

\medskip
\medskip
\medskip

{\bf 3. Spectrum of Hamiltonian of $(p,q;\alpha,
\beta,l)$-deformed oscillator}
\medskip
\medskip
\medskip

The Hamiltonian of the $(p,q;\alpha,\beta,l)$-deformed oscillator
system is defined in the same way as in case of the $q$-deformed
oscillator. From the relations
\begin{equation}
a a^+ - q^{-l} a^+ a = p^{\alpha N + \beta},\quad a a^+ - p^l a^+ a =
q^{-\alpha N - \beta}\end{equation} we have
\begin{equation} aa^+ |n\rangle = \frac{p^{-\alpha N-\beta
-l} - q^{\alpha N+\beta +l}}{p^{-l}-q^l}|n\rangle,\quad
a^+a|n\rangle = \frac{p^{-\alpha N-\beta }-q^{\alpha N+\beta
}}{p^{-l}-q^l}|n\rangle.\end{equation}

The Hamiltonian
\begin{equation}\label{burban: ham} H = a^+a + aa^+
\end{equation} of this $(p,q;\alpha,\beta,l)$-deformed
oscillator has a diagonal form in the basis $\{|n\rangle\}$:
\begin{equation}
H |n\rangle = \lambda_n|n\rangle,
\end{equation}
where
\begin{equation}\label{burban: spectr}\lambda_n =
\frac{p^{-\alpha n-\beta -l}-q^{\alpha n + \beta +l}}{p^{-l}- q^l}
+ \frac{p^{-\alpha n-\beta}-q^{\alpha n+\beta}}{p^{-l} - q^l}.
\end{equation}

Because of the identity $$\frac{p^{-\alpha n-\beta -l}- q^{\alpha
n + \beta +l}}{p^{-l}- q^l} = \frac{p^{-\alpha n -\beta -l}-
p^{-\alpha n-\beta }q^l + p^{-\alpha n-\beta }q^l - q^{\alpha n +
\beta + l}}{p^{-l} - q^l} $$
\begin{equation} = \frac{p^{-\alpha n-\beta }(p^{-l}- q^l) + (p^{-\alpha n-\beta }-
q^{\alpha n+\beta})q^l}{p^{-l}-q^l} = p^{-\alpha n - \beta} + q^l
\left({\frac{p^{-\alpha n-\beta}-q^{\alpha n + \beta}}{p^{-l} -
q^l} }\right)
\end{equation}
the relation (\ref{burban: spectr}) can be rewritten as
\begin{equation}
\label{burban: sym1} \lambda_n = p^{-\alpha n-\beta} + (q^l +
1)\left(\frac{p^{-\alpha n-\beta} - q^{\alpha n+\beta }}{p^{-l} -
q^l}\right).
\end{equation}

On the other hand $$ \frac{p^{-\alpha n-\beta - l}-q^{\alpha n +
\beta +l}}{p^{-1}-q^l} = \frac{p^{-\alpha n - \beta -l}-
p^{-l}q^{\alpha n + \beta} + p^{-l}q^{\alpha n + \beta} -
q^{\alpha n + \beta + l}}{p^{-1} - q^l}$$


\begin{equation}
= p^{-l}\left({\frac{p^{-\alpha n - \beta} - q^{\alpha n +
\beta}}{p^{-l}-q^l}}\right)+ q^{\alpha n + \beta},
\end{equation}
that is

\begin{equation}
\label{burban: sym2} \lambda_n = q^{\alpha n+\beta} + (p^{-l} +
1)\left(\frac{p^{-\alpha n-\beta} - q^{\alpha n+\beta }}{p^{-l} -
q^l}\right).
\end{equation}
It follows from (\ref{burban: sym1}) and  (\ref{burban: sym2})
that spectrum of the Hamiltonian (\ref{burban: ham}) is symmetric
under the change of parameter $q\to p^{-1},\quad p\to q^{-1}.$

\medskip
\medskip
\medskip
\medskip
{\bf 4. Hopf algebra structure of $(p,q;\alpha,\beta,l)$-deformed
oscillator algebra}
\medskip
\medskip
\medskip

It would be desirable to show that the generalized Heisenberg-Weyl
algebra, generated by the generators  defined above and the
relations (\ref{burban: basic1b}) and (\ref{burban: basic2b})
carries a Hopf algebra structure.

Remind, the associative algebra $C$ is a Hopf algebra if it admits
operations of homomorphisms of a coproduct $\Delta,$ a counit
$\epsilon $ and an anti-homomorphism of an antipode $S$:
\begin{equation} \Delta: C \to C\otimes \,C, \quad \Delta(a\, b) = \Delta(a)\Delta
(b),\end {equation} \begin{equation} \epsilon: C \to C, \quad
\epsilon(a b) = \epsilon(a)\,\epsilon(b)
\end{equation}
\begin{equation}
S(a\,b)= S(b)\,S(a).
\end{equation}
which satisfy properties
\begin{equation}
\label{burban: coalg1} (id \otimes \Delta)\Delta(h)=(\Delta
\otimes id) \Delta)\Delta(h),
\end{equation}
\begin{equation}\label{burban: coalg2}
(id\otimes\epsilon)\Delta(h)=(\epsilon \otimes id)\Delta(h),
\end{equation}
\begin{equation}\label{burban: coalg3}
m(id\otimes S)\Delta)(h)= m(S\otimes id)\Delta(h) =\epsilon(h){\bf
1}
\end{equation}
for all $h\in C.$

In our case the algebra is generated by $ {\bf 1}, a^+, a, N,$
satisfying  the relations
\begin{equation}\label{burban: alg1} [N, a] = -l a,\quad [N, a^+]
= l a^+,
\end{equation}
\begin{equation}
\label{burban: alg2} [a, a^+]_{A} = \frac{p^{-\alpha N - \beta_1
-l}- q^{\alpha N - \beta_1}}{p^{-l}- q^l}-{A}\frac{p^{-\alpha N -
\beta_2}- q^{\alpha N - \beta_2}}{p^{-l}- q^l}, \end{equation} and
a constant $A$ will be determined later on.

In particular, for $\beta_1 - \beta_2 = l$ we obtain the relation
(\ref{burban: rel3}) and at $p = q,\, l = 1$ this algebra reduced
to the one of \cite{OS}.

We define an action of coproduct $\Delta,$ counit $\epsilon,$ and
antipode $S$ on the generators of the algebra as
\begin{equation} \Delta(a^+) = c_1 a^+ \otimes p^{-\alpha_1 N} +
c_2 q^{\alpha_2 N }\otimes a^+, \end{equation}
\begin{equation}\Delta(a) = c_3 a\otimes p^{-\alpha_3 N} + c_4 q^{\alpha_4
N }\otimes a,\end{equation}
\begin{equation}\label{burban: number} \Delta( N ) = c_5 N \otimes {\bf 1} + c_6 {\bf 1}\otimes N +
\gamma{\bf 1 }\otimes{\bf 1},\end{equation}
\begin{equation}\Delta(\bf 1) = {\bf 1}\otimes{\bf 1},\end{equation}
\begin{equation}
\epsilon(a^+) = c_7,\quad \epsilon(a)= c_8, \end{equation}
\begin{equation} \epsilon (N)= c_9,\quad \epsilon ({\bf 1}) =
1,\end{equation}

\begin{equation}
S(a^+)= - c_{10} a^+,\quad S(a)= - c_{11}a,
\end{equation}
\begin{equation} S(N)= - c_{12}N + c_{13}{\bf 1},\quad S({\bf 1})= {\bf
1},
\end{equation}
where $c_i, i = 1\ldots 13$, and $\gamma$ are unknown coefficients
which must be determined by means of the rules of Hopf algebra
structure.

Using the relations
\begin{equation} a^+ \,r^{\alpha N} = r^{-\alpha l}\,r^{\alpha
N}a^+,\quad a \,r^{\alpha N} = r^{\alpha l}\, r^{\alpha N} a,
\end{equation} where $ r = p,q $
and
\begin{equation} q^{\alpha \Delta(N)} = q^{\alpha\gamma}q^{\alpha
N}\otimes q^{\alpha N}, \end{equation} we shall verify the axiom
(\ref{burban: coalg1}) of the Hopf algebra structure for $ h =
a^+, a, N.$ The condition (\ref {burban: coalg1}) for $h = a^+ $
gives $$(id\otimes\Delta)\Delta a^+ = c_1 p^{-\alpha_1\gamma}a^+
\otimes p^{-\alpha_1 N}\otimes p^{-\alpha_1 N} $$

\begin{equation}
\label{burban: first1} + c_1 c_2q^{\alpha_2 N}\otimes a^+\otimes
p^{-\alpha_1 N}+ c_2 c_2 q^{\alpha_2 N}\otimes q^{\alpha_2
N}\otimes a^+
\end{equation}
and
$$ (\Delta \otimes id) a^+  = c_1 c_1 a^+ \otimes p^{-\alpha_1
N}\otimes p^{\alpha_1 N} $$

\begin{equation}\label{burban: first2}
+ c_1 c_2q^{\alpha_2 N}\otimes a^+ \otimes p^{-\alpha_1 N} + c_2
q^{\alpha_2\gamma}q^{\alpha_2 N}\otimes q^{\alpha N}\otimes a^+.
\end{equation}
From (\ref{burban: first1}) and (\ref{burban: first2}) it follows
\begin{equation}
c_1 = p^{-\alpha_1\gamma},\quad c_2 = q^{\alpha_2\gamma}.
\end{equation}
The condition (\ref{burban: coalg1}) for $ h = a,$ $ N $ gives
\begin{equation} c_3 = p^{-\alpha_3\gamma},\quad c_4
= q^{\alpha_4\gamma}, c_5 = 1,\quad c_6 = 1.
\end{equation}

It is easy to see that
\begin{equation} \Delta(a)\Delta{(a^+)}= c_1 c_3
aa^+ \otimes p^{-(\alpha_1 + \alpha_3)N} + c_2 c_4 q^{(\alpha_2
+\alpha_4)N}\otimes a a^+ $$

$$+ c_2 c_3q^{\alpha_2l} q^{\alpha_2 N}a\otimes p^{-\alpha_3 N}a^+
+ c_1c_4p^{-\alpha_1 l}q^{\alpha_4 N} a^+ \otimes p^{-\alpha_1 N}
a\end{equation} and
\begin{equation} \Delta(a^+)\Delta (a) = c_1c_3 a^+a \otimes
p^{-(\alpha_1 + \alpha_3)N} + c_2 c_4 q^{(\alpha_2
+\alpha_4)N}\otimes a^+ a $$

$$+ c_2 c_3 p^{\alpha_3 l} q^{\alpha_2 N} a\otimes p^{-\alpha_3 N}
a^+ + c_1 c_4 q^{-\alpha_4 l}a^+ p^{-\alpha_4 N}a^+ \otimes
p^{-\alpha_1 N} a.\end{equation}

The action of the operation $\Delta$ on the left hand side of
(\ref{burban: alg2}) gives
\begin{equation}\label{burban: hom1}
\Delta(a)\Delta (a^+) - {A}\, \Delta(a^+)\Delta (a)= c_1
c_3[a,a+]_{A} \otimes p^{-(\alpha_1 +\alpha_3)N}+ c_2 c_4
q^{(\alpha_2 + \alpha_4)N}\otimes [a,a^+]_{A},\end{equation} if $$
q^{\alpha_2 l}-{A}\, p^{\alpha_3 l} = 0 \quad {\rm and }\quad
p^{-\alpha_1 l } - {A}\, q^{\alpha_4 l} = 0. $$ It leads to $$A =
p^{-\alpha_3 l }q^{\alpha_2 l},\quad A = p^{-\alpha_1
l}q^{\alpha_4 l}$$ or $\alpha_1 = \alpha_3,$ and $\alpha_2 =
\alpha_4.$

Using the relation (\ref{burban: alg2}) $$[ a, a^+]_{A} =
\frac{p^{-\alpha N - \beta_1}- q^{\alpha N + \beta_1}}{p^{-l}-
q^l}-{A}\,\frac{p^{-\alpha N - \beta_2}- q^{\alpha N +
\beta_2}}{p^{-l}- q^l} $$ $$= \frac{(p^{-\beta_1} - {A}\,
p^{-\beta_2}) p^{-\alpha N} - (q^{\beta_1}-
{A}\,q^{\beta_2})q^{\alpha N}}{p^{-l}- q^l} $$ one can represent
the expression (\ref{burban: hom1}) in the form $$ \Delta(a)\Delta
(a^+) - {A}\, \Delta(a^+)\Delta (a)$$ $$= c_1c_3
\frac{(p^{-\beta_1} - {A}\,p^{-\beta_2}) p^{-\alpha
N}-(q^{\beta_1}-{A}\, q^{\beta_2})q^{\alpha N}}{p^{-l}-q^l}
\otimes p^{-(\alpha_1 + \alpha_3)N}  $$ $$ + c_2 c_4 q^{(\alpha_2
+ \alpha_4)N}\otimes\frac{(p^{-\beta_1} - {A}\,p^{-\beta_2})
p^{-\alpha N}-(q^{\beta_1} - {A}\,q^{\beta_2})q^{\alpha
N}}{p^{-l}-q^l}$$ $$ = c_1c_3 \frac{(p^{-\beta_1} -
{A}\,p^{-\beta_2}) p^{-\alpha N}\otimes p^{-(\alpha_1 +
\alpha_3)N} - (q^{\beta_1}- {A}\,q^{\beta_2}) q^{\alpha N}\otimes
p^{-(\alpha_1 + \alpha_3)N}}{p^{-l}-q^l} $$
\begin{equation}\label{burban: hom2}+ c_2 c_4 \frac{(p^{-\beta_1}
-{A}\, p^{-\beta_2}) q^{(\alpha_2 + \alpha_4)N}\otimes p^{-\alpha
N}-(q^{\beta_1}- {A}\,q^{\beta_2})q^{(\alpha_2 + \alpha_4)
N}\otimes q^{\alpha N}}{p^{-l}-q^l}.\end{equation}

On the other hand, the action of the $\Delta$ on the right hand
side of (\ref{burban: alg2}) gives $$ \Delta(\frac {p^{-\alpha
N-\beta_1}-q^{\alpha N+\beta_1} - {A}\,( p^{-\alpha N-\beta_2}+
q^{\alpha N+\beta_2})}{p^{-l}-q^l}) $$
\begin{equation}\label{burban: hom3} = \frac{p^{-\alpha\gamma}( p^{-\beta_1}- {A}p^{-\beta_2})
p^{-\alpha N}\otimes p^{-\alpha N}-
q^{\alpha\gamma}(q^{\beta_1}-{A}\,q^{\beta_2})q^{\alpha N}\otimes
q^{\alpha N}}{p^{-l}-q^l}.\end{equation}

From (\ref{burban: hom2}) and (\ref{burban: hom3}) we have

$$ -c_1c_3(q^{\beta_1}- {A}\,q^{\beta_2})+
c_2c_4(p^{-\beta_1}-{A}\, p^{-\beta_2}) = 0. $$ If $\alpha_1 =
\alpha_2 = \alpha_3 = \alpha_4 = \alpha/2,$ then

\begin{equation} A = (p^{-1}q)^{\alpha l/2},\quad c_1c_3 = p^{-\alpha
\gamma},\quad \,\,c_1 c_4 = q^{\alpha\gamma}.\end{equation} It
follows $$ p^{-\alpha\gamma} (q^{\beta_1}-{A}\, q^{\beta_2})-
q^{\alpha\gamma}(p^{-\beta_1}-{A}\, p^{-\beta_2}) = 0, $$

\begin{equation}\label{burban: last}(pq)^{\alpha\gamma} =
\frac{q^{\beta_1}-{A}\,q^{\beta_2}}{p^{-\beta_1}-{A}\,p^{-\beta_2}}.
\end{equation}

The last equation (\ref{burban: last}) defines the parameter
$\gamma$ in the equation (\ref{burban: number}) for the Hopf
algebra structure.

Comparing the right-hand sides of the relations
\begin{equation}
\label{burban: coalg2'} (id\otimes \epsilon)\Delta(a^+) = c_1 a^+
\otimes p^{-\alpha_1 c_9} + c_2 q^{\alpha_2 N}\otimes c_7 +
\gamma{\bf 1}\otimes{\bf 1},\end{equation} and
\begin{equation}\label{burban: coalg"}
(\epsilon \otimes id )\Delta(a^+) = c_1 c_7 {\bf 1}\otimes
p^{-\alpha_1 N} + c_2 q^{\alpha_2 c_9}{\bf 1}\otimes a^+ + \gamma
{\bf 1}\otimes {\bf 1}\end{equation} and using the axiom
(\ref{burban: coalg2}) for the generator $a^+,$ we obtain (take
into account that $c_1 = p^{-\alpha_1 \gamma},\quad c_2 =
q^{\alpha_2 \gamma}$) $$c_1 p^{-\alpha_1 c_9} = c_2 q^{\alpha_2
c_9},\quad - \alpha_1\gamma - \alpha_1 c_9 = 0,$$ hence
\begin{equation}  c_9 = - \gamma.\end{equation}

An easy calculation gives
\begin{equation}\label{burban: last1} m(id\otimes S)\Delta
(a^+)= c_1 p^{-\alpha_1 c_{13}} a^+ p^{\alpha_1 c_{12} N} - c_2
c_{10}q^{\alpha_2 N} a^+
\end{equation}
and
\begin{equation}\label{burban: last2}
m(S\otimes id)\Delta (a^+) = - c_1 c_{10} a^+ p^{-\alpha_1 N} +
c_2 q^{\alpha_2 c_{13}} q^{-\alpha_2 c_{13} N} a^+.
\end{equation}
From these relations and from the axiom (\ref{burban: coalg3}) for
$a^+$ we obtain $$ p^{-\alpha_1c_{13}}a^+ p^{\alpha_1c_{12} N} = -
c_{10} a^+ p^{-\alpha_1 N},\quad - c_{10}q^{\alpha_2 N} =
q^{\alpha_2 N}q^{-\alpha_2 c_{13}N}. $$ Then
\begin{equation} c_{10}= -1,\,c_{12}= -1,\, c_{13} =
0.\end{equation}

The same calculations for $a$ give $c_{11}= -1.$

A fulfillment of the remainder relations of the algebra under the
action of the Hopf algebra operations can be easily verified.

\medskip
\medskip
\medskip

\subsection*{Acknowledgements}
I would like to thank A.U. Klimyk for many useful discussions,
valuable suggestions.

\noindent This research was partially supported by Grant 10.01/015
of the State Foundation of Fundamental Research of Ukraine.


\medskip
\medskip
\medskip



\begin{thebibliography}{99}

\bibitem {ACL} Arik, D. D., Coon Y., and Lam A.,
{\it J. Math. Phys.,} 1975, V 16, 1776.

\bibitem{CZ} Curtright T., Zachos C., {\it Phys. Lett. B,} 1990, V 243, 237.

\bibitem {ShWZ} Schirmacher, J., Wess J., Zumino B., {\it Z. Phys. C,} 1992,
V 49, 317.

\bibitem {SW} Smirnov Yu. F., Werhan R.F., {\it J.Phys. A,} 1992, V 25, 5563.

\bibitem {B} Biedenharn L. C., The spectrum group $SU_q(2)$ and a q-analogue of the boson
operator, {\it J. Phys. A: Math. Gen.}, 1989, V 22, L873 - L878.

\bibitem {M} Macfarlane A. J., On q-analogue of the quantum harmonic oscillator and
quantum droup $SU_q(2)$, {\it J. Phys. A,} 1989, V22, 4581 - 4585.

\bibitem{ChJ} R.Chacrabarti, Jagannatan, A $(p,q)$-oscillator realization of
two-parameter quantum algebras, {\it J. Phys A: Math. Gen.,} 1991,
V24, L711 - L718.

\bibitem {J} Jing S., {\it Nuovo Cimento A,} 1995, V 105, 1267.

\bibitem {Q} Quesne C., Two-parameter versus of one-parameter quantum deformation of
$su(2)$, {\it Phys. Lett. A}, 1993, V174, 19 - 24.

\bibitem {BK} Burban I. M., Klimyk A. U.,
P,Q-differentiation, P,Q-integration, and P,Q-hypergeometric
functions related to quantum groups, {\it Integral Transforms and
Special Functions}, 1994, V2, 15 - 36.

\bibitem {JR} Jaghannatan R., Rao K.S., Two-parameter quantum algebras, twin-basic numbers,
and associated hypergeometric series, {\it arXiv math.
NT/0602613}.

\bibitem{Br} Burban I. M., Two-parameterdeformation of oscillator algebra, {\it Phys. Lett. B}
 1993, V 319, 485 - 489.

\bibitem{C-C-N-U} Chung W., Chung K., Nam S-T., Um C., Generalized deformed algebra, {\it Phys. Lett.
A} 1993, V 183, 363-370.

\bibitem{OS} Oh C.H., Sing K., Generalizedc q-oscillators and thier
Hopf structures, {\it J. Phys. A: Math. Gen}, 1994, V 27, 5907 -
5918.

\bibitem {BK1} Burban I. M., Klimyk A. U., On spectral properties of $q$- oscillator operators, {\it Lett. Math. Phys.} 1991, V
29, 13 - 18.

\bibitem {BDY} Borzov V. V., Damaskinsky E. V., Yagorov S.V., Some
representations of the generalized deformed oscillator algebra,
{\it arXiv q-alg/9509022}.

\end{thebibliography}
\end{document}